
\documentstyle{amsppt}
\magnification=\magstep1
\NoBlackBoxes
\settabs 10 \columns

\def\qed {\ \ \ \vrule height6pt  width6pt depth0pt}
\def\L{\Cal L}
\def\N{\Bbb N}
\def\C{\Bbb C}
\def\P{\Cal P}
\def\F{\Cal F}
\def\A{\Cal A}
\def\S{\Cal S}
\def\M{\Cal M}
\def\H{\Cal H}
\def\U{\Cal U}
\def\V{\Cal V}
\def\K{\Cal K}
\def\G{\Cal G}
\def\Nu{\Cal N}
\topmatter

\title Factorization and Reflexivity on Fock spaces
\endtitle

\author Alvaro Arias and Gelu Popescu
\endauthor

\thanks
The first author was supported in part by NSF DMS 93-21369
\endthanks

\address A. Arias
\newline Division of Mathematics, Computer Science and Statistics,
The University of Texas at San Antonio, San
Antonio, TX 78249, U.S.A.\endaddress
\email arias\@ringer.cs.utsa.edu\endemail

\address G. Popescu
\newline Division of Mathematics, Computer Science and Statistics,
The University of Texas at San Antonio, San
Antonio, TX 78249, U.S.A.\endaddress
\email gpopescu\@ringer.cs.utsa.edu\endemail

\abstract
The framework of the paper is that of the full Fock space
${\Cal F}^2({\Cal H}_n)$ and the Banach algebra $F^\infty$ which can be
viewed
as non-commutative analogues of the Hardy spaces $H^2$ and
$H^\infty$ respectively.

An inner-outer factorization for any element in
${\Cal F}^2({\Cal H}_n)$ as
well as characterization of invertible elements in $F^\infty$ are
obtained. We also give a complete characterization of invariant
subspaces for the left creation operators $S_1,\cdots, S_n$ of
${\Cal F}^2({\Cal H}_n)$.
This enables us to show that every weakly (strongly)
closed unital subalgebra of $\{\varphi(S_1,\cdots,S_n):\varphi\in
F^\infty\}$
is reflexive, extending in this way the classical result of Sarason
[S].

Some properties of inner and outer functions and many examples are also
considered.
\endabstract


\endtopmatter

\document

\baselineskip18pt

\heading
0. Introduction
\endheading

Let $H^p$, $1\leq p\leq\infty$ be the classical Hardy spaces on the
disk.
$H^2$ is a Hilbert space with orthonormal basis $\{z^n\}_{n=0}^\infty$
and
it is well known that if $\varphi\in H^\infty$, then
$\|\varphi\|_\infty=\|M_\varphi\|$, where $M_\varphi:H^2\to H^2$ is the
multiplication operator. Hence,
$$\|\varphi\|_\infty=\sup~\{\|\varphi p\|_2:\|p\|_2\leq 1,
\hbox{ and $p$ is a polynomial on $z$}\}.\leqno{(1)}$$

The set of polynomials on $z$, $\P(z)$, determine the Hardy spaces:
$H^2$ is the closure of $\P(z)$ in the Hilbert space with orthonormal
basis
$\{z^n\}_{n=0}^\infty$. Once we have $H^2$, $H^\infty$ consists of all
$\varphi\in H^2$ such that (1) is finite.

\medbreak

In [Po3], [Po4] the second author introduced and studied a
``non-commutative''
analogue of $H^\infty$. Let $\P$ be the set of polynomials in
$n$ non-commutative indeterminates $e_1,\cdots, e_n$. To stress the
non-commutativity of the product we use tensor notation; (i.e.,
$e_1\otimes e_2$ represents the product of $e_1$ and $e_2$ and is
different from $e_2\otimes e_1$). A typical element of $\P$ looks
like
$$p=a_0+\sum_{k=1}^m\sum_{1\leq i_1,\cdots, i_k\leq n}
a_{i_1\cdots i_k}e_{i_1}\otimes\cdots\otimes e_{i_k},
\leqno{(2)}$$
where $m\in\N$ and $a_0, a_{i_1\cdots i_k}\in\C$. Let $\F^2(\H_n)$ be
the
Hilbert space having the monomials as an orthonormal basis (i.e., $1$
and the
elements of the form $e_{i_1}\otimes\cdots\otimes e_{i_k}$).
This will be viewed as analogue to $H^2$;   $\P$ is dense
there.

Define $F^\infty$ as the set of all $g\in\F^2(\H_n)$ such that
$$\|g\|_\infty=\sup~\{\|g\otimes p\|_2: p\in\P, \|p\|_2\leq1\}<\infty,
\leqno{(3)}$$
where $\|\cdot\|_2=\|\cdot\|_{\F^2(\H_n)}$.

It is easy to prove (see [Po3]) that $(F^\infty,\|\cdot\|_\infty)$ is a
non-commutative Banach algebra (if $f,g\in F^\infty$, then
$f\otimes g\in F^\infty$ and
$\|f\otimes g\|_\infty\leq\|f\|_\infty\|g\|_\infty$)
which can be viewed as a non-commutative analogue of the Hardy space
$H^\infty$; when $n=1$ they coincide.

\medbreak

$F^\infty$ shares many properties with $H^\infty$.
For instance, the second author (see [Po3])
extended the classical von Neumann's inequality [vN]
to polynomials in $F^\infty$. He proved that if $T_1,\cdots, T_n\in
B(\H)$,
for some Hilbert space $\H$, are such that the operator matrix
$[T_1\cdots T_n]$ is a contraction (i.e., $\|\sum_{i=1}^n T_i
T_i^*\|\leq1$)
then for every $p\in\P$,
$$\|p(T_1,\cdots,T_n)\|\leq\|p\|_\infty,$$
where $p(T_1,\cdots,T_n)$ is obtained from (2) by replacing the $e_i$'s
by
$T_i$'s.

There are also many analogies with the invariant subspaces of
$H^2$, inner and
outer functions in $H^\infty$,
Toeplitz operators, etc.  (see [Po2], [Po4]). We
will give precise definitions for many of these facts below.

\medbreak

For $i\leq n$ define $S_i:\F^2(\H_n)\to\F^2(\H_n)$ by
$S_i\psi=e_i\otimes\psi.$
These are unilateral shifts with orthogonal final spaces and
are analogue to $S=M_z:H^2\to H^2$, multiplication by $z$, the
unilateral shift on $H^2$.

For $\varphi\in F^\infty$ let
$\varphi(S_1,\cdots,S_n):\F^2(\H_n)\to \F^2(\H_n)$ be
defined by $\varphi(S_1,\cdots, S_n)\psi=\varphi\otimes \psi$ for
every $\psi\in \F^2(\H_n)$. Following the classical case the second
author
(see [Po4]) defined:

\item{1.} $\varphi\in F^\infty$ is {\sl inner} iff
$\varphi(S_1,\cdots,S_n)$ is
an isometry on $\F^2(\H_n)$, and

\item{2.} $\varphi\in F^\infty$ is {\sl outer} iff
$\varphi\otimes\F^2(\H_n)$ is dense in $\F^2(\H_n)$.

In this paper we define

\item{3.} $\psi\in \F^2(\H_n)$ is {\sl outer} iff
$\psi\otimes\P$ is dense in $\F^2(\H_n)$.

\medbreak
In Section 1 we will set the notation, terminology and prove a few
elementary facts about $\F^2(\H_n)$, $F^\infty$ and inner and
outer functions.
In Section 2 we will study some factorization
properties of inner and outer functions
and will complete the characterization of the invariant
subspaces for $S_1,\cdots, S_n$ started in [Po1].
In Section 3 we will provide many examples
of inner functions. We will see that the theory of inner functions
is more interesting in the non-commutative case than in the commutative
one. In Section 4 we will prove that every weakly (strongly) closed
unital subalgebra $\A$ of $F^\infty$ is reflexive; i.e., Alg
Lat$\,\A=\A$.
The proof of this uses most of the examples of inner functions
presented in Section 3. We will finish in Section 5 with some open
problems.

\heading
1. Preliminaries
\endheading

In this section we will set the notation, terminology and
prove a few elementary facts.

$\F^2(\H_n)$ is called the full Fock space on the $n$-dimensional
Hilbert
space $\H_n$ with orthonormal basis $(e_1,\cdots, e_n)$. It is usually
represented by
$$\F^2(\H_n)=\C1\oplus \bigoplus_{m\geq1}\H_n^{\otimes m}.$$

Let $\Lambda=\{1,2,\cdots, n\}$ be fixed
throughout the paper. For every $k\geq 1$, let $F(k,\Lambda)$ be the
set of
all functions from $\{1,2,\cdots, k\}$ to $\Lambda$, and let
$$\F=\bigcup_{k=0}^\infty F(k,\Lambda),$$
where $F(0,\Lambda)$ stands for $\{0\}$. We use $\F$ to describe
$\F^2(\H_n)$ and simplify the notation.
If $f\in F(k,\Lambda)$, let
$$e_f=e_{f(1)}\otimes e_{f(2)}\otimes\cdots\otimes e_{f(k)},
\quad\text{ and for }\quad k=0,\quad
e_0=1.$$
Then $\F^2(\H_n)$ is the Hilbert space with orthonormal basis
$\{e_f:f\in\F\}$.

Moreover, we use $\F$ to describe finite products.
If $T_1, T_2, \cdots, T_n\in B(\H)$ for some Hilbert space $\H$, and
$f\in F(k,\Lambda)$ we denote
$$T_f=T_{f(1)} T_{f(2)}\cdots T_{f(k)}, \quad\quad T_0=I.\leqno{(4)}$$

\medbreak

We will use the following elementary results.

\proclaim{Lemma 1.1} Let $\varphi, \psi\in F^\infty$. If
$\varphi\otimes \psi=0$, then $\varphi=0$ or $\psi=0$.
\endproclaim

\demo{\bf Proof}
Let $\varphi=\sum_{f\in\F}a_fe_f$ and
$\psi=\sum_{g\in\F}b_ge_g$ be non-zero elements.
Let $k$ is the smallest integer
such that for some $\bar{f}\in F(k,\Lambda)$, $a_{\bar{f}}\not=0$, and
$\ell$ the smallest integer such that for some $\bar{g}\in
F(\ell,\Lambda)$,
$b_{\bar{g}}\not=0$.

We have that $\varphi\otimes\psi=\sum_{h\in\F}
\bigl[\sum_{e_f\otimes e_g=e_h}a_fb_g\bigr]e_h$.
Let $\bar{h}$ be such that $e_{\bar{h}}=e_{\bar{f}}\otimes
e_{\bar{g}}$. It is easy to see that
$$\langle \varphi\otimes\psi,e_{\bar{h}}\rangle=
\sum_{e_f\otimes e_g=e_{\bar{h}}}a_fb_g=
a_{\bar{f}}b_{\bar{g}}\not=0.$$
Hence, $\varphi\otimes\psi\not=0$.\qed
\enddemo

\proclaim{Corollary 1.2} Suppose that $\varphi\in F^\infty$,
$\varphi\not=0$ and
that $\varphi=\varphi\otimes\psi$. Then $\psi=1$.
\endproclaim

\demo{\bf Proof}
If $\varphi=\varphi\otimes\psi$, then
$\varphi\otimes(1-\psi)=0$. Since $\varphi\not=0$, we get that
$1-\psi=0$.\qed
\enddemo

\medbreak

Recall that $\varphi\in F^\infty$ is inner iff the map
$\psi\to\varphi\otimes\psi$ is an isometry on $\F^2(\H_n)$ and that
$\psi\in F^\infty$ is outer iff $\psi\otimes\F^2(\H_n)$ is dense
in $\F^2(\H_n)$.
It is immediate from the definitions that

\proclaim{Proposition 1.3} $\varphi$ is inner iff $\{\varphi\otimes
e_f: f\in\F\}$
is an orthonormal set in $\F^2(\H_n)$.
\endproclaim

\proclaim{Proposition 1.4} $\psi$ is outer iff there exist $h_n\in
\F^2(\H_n)~$
(or $p_n\in\P$) such that $\psi\otimes h_n\to e_0$
($\psi\otimes p_n\to e_0$) in $\F^2(\H_n)$.
\endproclaim

We will also use

\proclaim{Proposition 1.5} Let $\varphi$ be inner. The map
$\psi\to\varphi\otimes\psi$ is an isometry on $F^\infty$. Moreover,
if $\phi=\varphi\otimes\psi$, then $\psi\in F^\infty$ if and only if
$\phi\in F^\infty$ and $\|\psi\|_\infty=\|\phi\|_\infty$.
\endproclaim

\demo{\bf Proof}
Let $\varphi$ be inner and $\psi\in F^\infty$.
For every $p\in\P$ one has
$$\|(\varphi\otimes\psi)\otimes p\|_2=\|\varphi\otimes(\psi\otimes
p)\|_2
=\|\psi\otimes p\|_2.$$
It follows from (3) that
$\|\varphi\otimes\psi\|_\infty=\|\psi\|_\infty.$

Suppose now that $\phi=\varphi\otimes\psi$. It is clear that
$\sup_{p\in(\P)_1}\|\phi\otimes p\|_2$ is finite if and only if
$\sup_{p\in(\P)_1}\|\psi\otimes p\|_2$ is finite too. \qed
\enddemo

\proclaim{Proposition 1.6} $\varphi\in F^\infty$ is inner iff
$\|\varphi\|_2=\|\varphi\|_\infty=1$.
\endproclaim

\demo{\bf Proof}
Suppose that $\varphi\in F^\infty$ is inner. Then
$\|\varphi\|_\infty\geq \|\varphi\otimes e_0\|_2=\|\varphi\|_2=1$. On
the
other hand, if $p=\sum_f a_f e_f\in\P$, we have that
$\|\varphi\otimes p\|_2^2=\sum_f|a_f|^2\|\varphi\otimes e_f\|_2^2=
\sum_f|a_f|^2=\|p\|_2^2.$
Therefore $\|\varphi\|_\infty\leq 1.$

Conversely, suppose that $\|\varphi\|_\infty=\|\varphi\|_2=1$. Then
for every $f\in\F$, $1=\|\varphi\|_2=\|\varphi\otimes e_f\|_2$.
It is well known that if $x,y\in\ell_2$ satisfy $\|x\|_2=\|y\|_2=1$ and
$\|x\pm y\|_2\leq\sqrt{2}$, then $x$ and $y$ are orthogonal. Let
$f,g\in\F$,
$f\not= g$. Since $\|\varphi\|_\infty=1$,
$\|\varphi\otimes e_f\pm \varphi\otimes e_g\|_2=
\|\varphi\otimes(e_f\pm e_g)\|_2\leq\|e_f\pm e_g\|_2=\sqrt{2}$.
Therefore,
$\varphi\otimes e_f\perp \varphi\otimes e_g$.\qed
\enddemo
\medbreak

A sequence $\S=\{S_\lambda\}_{\lambda\in\Lambda}$ of unilateral shifts
on a Hilbert space $\H$ with orthogonal final spaces is called a
$\Lambda$-orthogonal shift if the operator matrix
$[S_1 S_2 \cdots S_n]$ is non-unitary; i.e.,
$\L=\H\ominus\bigl(\bigoplus_{\lambda\in\Lambda}S_\lambda
\H\bigr)\not=\{0\}$.
$\L$ is called the {\sl wandering subspace} of $\S$. The second author
(see [Po2]) used $\L$ to get a Wold decomposition for $\H$. He proved
that
$$f,g\in\F,~ f\not=g\Longrightarrow S_f\L\perp S_g\L,\leqno{(5)}$$
(see (4) for the notation of $S_f$), and
$$\H=\bigoplus_{f\in\F}~ S_f\L.\leqno{(6)}$$
Moreover the dimension of $\L$, called the multiplicity of the shift,
determines up to unitary equivalence the $\Lambda$-orthogonal shift.

A model for the $\Lambda$-orthogonal shifts of multiplicity one
is provided by the unilateral shifts $S_i:\F^2(\H_n)\to\F^2(\H_n)$
defined by
$S_i\psi= e_i\otimes\psi$, $i\leq n$.
These operators are also called the left creation operators on the full
Fock space.
\medbreak

Let $U:\F^2(\H_n)\to \F^2(\H_n)$ be the flip operator defined by
$U\varphi=\tilde{\varphi}$ where
$$e_f=e_{f(1)}\otimes \cdots\otimes e_{f(k)},\quad\quad
Ue_f=\tilde{e}_f=e_{f(k)}\otimes e_{f(k-1)}\otimes\cdots\otimes
e_{f(1)}.$$
It is easy to check that $U$ is unitary and that
$U=U^*=U^{-1}$. We use $U$ to describe
right multiplication on $\F^2(\H_n)$.

\proclaim{Lemma 1.7} Let $\psi\in\F^2(\H_n)$ and suppose that
$$\sup\{\|p\otimes\psi\|_2~:~p\in\P,~\|p\|_2\leq1\}<\infty.$$
Then $\tilde{\psi}\in F^\infty$ and the supremum is equal to
$\|\tilde{\psi}\|_\infty.$
\endproclaim

\demo{\bf Proof}
Since $U$ is an isometry on $\F^2(\H_n)$ and
$U(p\otimes q)=U(q)\otimes U(p)=\tilde{q}\otimes\tilde{p}$, we have
$$\sup_{p\in(\P)_1}\|p\otimes\psi\|_2
=\sup_{p\in(\P)_1}\|\tilde{\psi}\otimes\tilde{p}\|_2
=\sup_{p\in(\P)_1}\|\tilde{\psi}\otimes p\|_2
=\|\tilde{\psi}\|_\infty.\qed$$
\enddemo

\noindent{\bf Remark.} $U$ is unbounded in $F^\infty$. Example 6 of
Section 3
tells us that if $p(e_1)$ is a polynomial in $e_1$, then
$\|e_2\otimes p(e_1)\|_\infty=\|p(e_1)\|_\infty$ but
$\|U(e_2\otimes p(e_1))\|_\infty=\|p(e_1)\otimes
e_2\|_\infty=\|p(e_1)\|_2$.

\heading
2. Factorization results
\endheading

In this section we will prove

\proclaim{Theorem 2.1}
If $~\psi\in F^2(H_n), \ \psi\not= 0$, then there exist $\varphi$ inner
and $g$ outer functions such that $~\psi=\varphi\otimes g$.
Moreover, the factorization is essentially unique and
$\psi\in F^\infty$ if and only if $g\in F^\infty$ and
$\|\psi\|_\infty=\|g\|_\infty$.
\endproclaim

And together with [Po3, Corollary 3.5],

\proclaim{Theorem 2.2} Let $~\Phi:F^\infty\to B(\F^2(H_n))$ be defined
by
$\Phi(\varphi)=\varphi(S_1,\cdots,S_n)$. Then:
\roster
\item"{(i)}" $\Phi\ $ is an algebra homomorphism.
\item"{(ii)}" $\|\Phi \psi\|=\|\psi\|_\infty$, for any $~\psi\in
F^\infty$.
\item"{(iii)}" $\psi~$ is invertible in $~F^\infty$
if and only if $~\Phi \psi~$ is invertible in $~B(\F^2(H_n))$.
\endroster
\endproclaim

To prove Theorem 2.1 we will first complete the characterization of
invariant subspaces for   $S_1,\cdots,S_n$ started in [Po1]
for $\Lambda$-orthogonal shifts of arbitrary multiplicity.
\smallbreak

We say that the inner functions $\varphi_1$, $\varphi_2$ are
{\sl orthogonal} if $\varphi_1\otimes\F^2(\H_n)\perp
\varphi_2\otimes\F^2(\H_n)$.
\proclaim{Theorem 2.3} If $\M\subset\F^2(H_n)$ is invariant for each
$S_1,\cdots,S_n$ then there exists a sequence $\{\varphi\}_{j\in J}$
of orthogonal inner functions
such that
$$\M=\bigoplus_{j\in J}F^2(H_n)\otimes\tilde{\varphi}_j.$$
Moreover, this representation is essentially unique.
\endproclaim
\demo{\bf Proof}
Let $\M$ be a nontrivial closed invariant subspace for
$S_1,\cdots,S_n$.
According to (6)
$$\M=\bigoplus_{f\in\F}S_f\L,\leqno{(7)}$$
where $\L=\M\ominus(S_1\M\oplus\cdots\oplus S_n\M)$ is the wandering
subspace
for $S_1|_{\M},\cdots,S_n|_{\M}$.

Let $\{\tilde{\varphi}_j\}_{j\in J}$ be an orthonormal basis in $\L$.
According  to (5), for each $j\in J$
$$ S_f\tilde{\varphi}_j \perp S_g\tilde{\varphi}_j,\quad
\text{ for any}\quad f,g\in \F, f\not=g.\leqno{(8)}$$
By Proposition 1.3 we infer that $\varphi_j$ is an inner function.
On the other hand the relations (7) and (8) imply the
following orthogonal decomposition
$$\M=\bigoplus_{j\in J}F^2(H_n)\otimes\tilde{\varphi}_j.$$

Suppose that $\M=\bigoplus_{i\in I}[\F^2(\H_n)\otimes\tilde{\psi}_i]$
for some orthogonal inner functions $\psi_i$. Let
$\L'=\vee_{i\in I}\tilde{\psi}_i$. It is easy to see that
$$\M=\bigoplus_{f\in S_f}~S_f\L'.$$
This implies that $\L'=\M\ominus(S_1\M\oplus\cdots\oplus S_n\M)$ is
the wandering subspace for $S_1|_{\M},\cdots,S_n|_{\M}$ and
$\L=\L'$. Since $\{\tilde{\varphi}_j\}_{j\in J}$ and
$\{\tilde{\psi}_i\}_{i\in I}$ are orthogonormal basis
in $\L$ we infer that they have the same cardinality and that there
exists a unitary operator $U\in B(\L)$ such that
$U\tilde{\varphi}_j=\tilde{\psi}_j$, $j\in J$. \qed

\enddemo
\medbreak

\proclaim{Corollary 2.4} If $\varphi_1,~\varphi_2$ are inner functions
such
that $\varphi_1\otimes\F^2(\H_n)=\varphi_2\otimes\F^2(\H_n)$, then
there exists $\alpha\in\C$, $|\alpha|=1$ such that
$\varphi_1=\alpha\varphi_2$.
\endproclaim

\demo{\bf Proof} Let $\M=\F^2(\H_n)\otimes\tilde{\varphi}_1=
\F^2(\H_n)\otimes\tilde{\varphi}_2$. It is clear that $\M$ is
invariant for $S_1,\cdots, S_n$. It follows from the proof
of Theorem 2.3 that $[\tilde{\varphi}_1]=[\tilde{\varphi}_2]$.
Since $\tilde{\varphi}_1$, $\tilde{\varphi}_2$
are normalized in $\F^2(\H_n)$ there exists $\alpha\in\C$, $|\alpha|=1$
such that $\tilde{\varphi}_1=\alpha \tilde{\varphi}_2$.\qed
\enddemo

The next lemma will show that all the cyclic invariant subspaces
of $\F^2(\H_n)$ are of the form $\F^2(\H_n)\otimes\tilde{\varphi}$ for
some inner function $\varphi$.

\proclaim{Lemma 2.5}
Let $~\V=\{V_1,\cdots,V_n\}~$ be a $\Lambda$-orthogonal shift of
arbitrary multiplicity acting on the Hilbert space $~\K$.
If $~\M\subset\K~$ is a cyclic invariant subspace of $~\V$, that is
$$
\M=\bigvee_{f\in\F} V_f \psi,\quad\text{for some\ }
\psi\in\K,\ \psi\not=0
$$
then $\V_{|\M}:=\{V_1{|_\M},\cdots,V_n{|_\M}\}$
is unitarily equivalent to the $\Lambda$-orthogonal shift of
multiplicity 1.
\endproclaim

\demo{\bf Proof}
Let $~\psi\in\K,\  \psi\not= 0~$ be a fixed element in $~\K~$ and let
$$
\M=\bigvee_{f\in\F}V_f \psi.
$$
Since $~\V~$ is a $\Lambda$-orthogonal shift in $~\K~$ it follows from
(5) and (6) that
$\K=\bigoplus_{f\in\F}V_f\L,$
where $\L=\K\ominus(V_1\K\oplus\cdots\oplus\V_k\K)$
is the wandering subspace of $~\V$.

Since $~\V|_{\M}~$ is also a $\Lambda$-orthogonal shift let us denote
by
$~\L_0~$ its wandering subspace.
Therefore, $~\L_0\subset\M~$ and
$$
\M=\bigoplus_{f\in\F}~V_f\L_0,\quad\text{where}
\quad\L_0=\M\ominus(V_1\M\oplus\cdots\oplus V_n\M).\leqno{(9)}
$$

Consider $~w_0=P_{\L_0}\psi$, where $~P_{\L_0}~$ is the orthogonal
projection from $~\K~$ onto $~\L_0$.
Notice that according to (9) it follows that $~w_0\not=0$.
Let $~\ell_0\in\L_0~$ be such that $~\ell_0\perp w_0$.
Since $~\L_0~$ is wandering subspace for $~\V{|_\M}~$ we have
$$
\ell_0\perp V_f \psi,\quad\text{for any\ }f\in F(k,\Lambda),\ k\ge 1.
$$
On the other hand
$$
\langle \ell_0, \psi\rangle= \langle \ell_0,\ P_{\L_0} \psi\rangle=
\langle \ell_0, w_0\rangle=0.
$$
Therefore, $~\ell_0\perp V_f \psi,\ \text{\ for any\ } f\in\F$.
Since $~\ell_0\in\M=\bigvee\limits_{f\in\F} V_f \psi~$
we infer that $~\ell_0=0$.
Thus, $\ \text{dim}~\L_0=1~$ and according to [Po3, Theorem 1.2]
it follows that $~\V|_{\M}$ is unitarily equivalent to
the $~\Lambda$-orthogonal
shift of multiplicity one
$~\S=\{S_1,\cdots,S_n\}$.
The proof is complete. \qed
\enddemo

\proclaim{Corollary 2.6}
If $\psi\in F^2(H_n),~\psi\not=0$ then there exists an inner function
$\varphi \in F^\infty~$ such that
$$\text{clos}[\P\otimes \psi]=F^2(H_n)\otimes\tilde{\varphi}.$$
Moreover, the representation is essentially unique.
\endproclaim

\demo{\bf Proof} Let $\M=\text{clos}[\P\otimes \psi]$. It follows from
the previous lemma that $\L=\M\ominus(S_1\M\oplus\cdots\oplus S_n\M)$,
the wandering subspace for $S_1|_{\M},\cdots,S_n|_{\M}$, has dimension
1.
Using Proposition 2.3 we get that $\M=\F^2(\H_n)\otimes\tilde{\varphi}$
for some inner function $\varphi$. The ``uniqueness'' is a consequence
of Corollary 2.4.\qed
\enddemo

\demo{\bf Proof of Theorem 2.1}
Let $~\M=\text{clos~}[\P\otimes\tilde{\psi}]~$.
Then $~\M~$ is a nontrivial
closed invariant subspace for each $~S_1,\cdots,S_n~$.
By Corollary 2.6  it is of the form
$\M=F^2(H_n)\otimes\tilde{\varphi}$, where $~\varphi\in F^\infty~$ is
inner.
Since $~\tilde{\psi}\in\M$, there must exist $~\tilde{g}\in F^2(H_n)~$
such that
$$
\tilde{\psi}=\tilde{g}\otimes\tilde{\varphi}
\quad\text{ and then }\quad \psi=\varphi\otimes g.
$$

Since $\tilde{\varphi}\in\M$, there exist $p_n\in\P$ such that
$$\tilde{\varphi}=\lim_{n\to\infty} p_n\otimes\tilde{\psi}
=\lim_{n\to\infty}(p_n\otimes\tilde{g})\otimes\tilde{\varphi}
=\bigl[\lim_{n\to\infty}p_n\otimes\tilde{g}\bigr]\otimes\tilde{\varphi}.$$
We use that $\varphi$ is inner for the last equality. By Corollary 1.2,
$\lim_{n\to\infty}p_n\otimes\tilde{g}=e_0$. Hence
$\lim_{n\to\infty}g\otimes\tilde{p}_n=e_0$, and using Proposition 1.4
we conclude that $g$ is outer. It follows
from Proposition 1.5 that $g\in F^\infty$ iff $\psi\in F^\infty$.

Suppose now that $\psi=\varphi_1\otimes g_1=\varphi_2\otimes g_2$
for some $\varphi_1,~\varphi_2$ inner and $g_1,~g_2$ outer. Then
we have that $\varphi_1\otimes\F^2(\H_n)=\varphi_2\otimes\F^2(\H_n)$,
and
using Corollary 2.4 we get that $\varphi_1=\alpha\varphi_2$ for some
$\alpha\in\C$, $|\alpha|=1$. Hence
$0=\varphi_1\otimes g_1-\varphi_2\otimes g_2
=\varphi_2\otimes(\alpha g_1-g_2)$. By Lemma 1.1,
$\alpha g_1=g_2$.\qed
\enddemo

As in the classical case, there is a factorization result for inner
functions.

\proclaim{Corollary 2.7} Let $\varphi_1,~\varphi_2$ be inner
functions.
$\varphi_1\otimes\F^2(\H_n)\subset\varphi_2\otimes\F^2(\H_n)$ if
and only if there exists $\varphi_3$ inner such that
$\varphi_1=\varphi_2\otimes\varphi_3$.
\endproclaim

\demo{\bf Proof} Suppose that
$\varphi_1\otimes\F^2(\H_n)\subset\varphi_2\otimes\F^2(\H_n)$.
Then $\varphi_1=\varphi_2\otimes\psi$ for some $\psi\in\F^2(\H_n)$.
Let $\psi=\psi_i\otimes\psi_e$ be the inner-outer factorization of
$\psi$.
Therefore, $\varphi_1=(\varphi_2\otimes\psi_i)\otimes\psi_e$.
By the uniqueness of the factorization of $\varphi_1$ we deduce that
$\psi_e=1$. The converse is clear.\qed
\enddemo

\bigbreak

We will finish the details of Theorem 2.2 now.
Parts (i) and (ii) appear in [Po3]. The third part
is included in the proof of the next theorem.

\proclaim{Theorem 2.8}
$\varphi\in F^\infty$ is invertible if and only if
$\varphi$ is outer and there exists
$\delta>0$ such that
$\|\varphi\otimes p\|_2\ge\delta\|p\|_2$, for any $p\in\P.$
\endproclaim

\demo{\bf Proof.}
Let $\varphi\in F^\infty$ be invertible and find
$\psi\in F^\infty$ such that
$\varphi\otimes\psi=\psi\otimes\varphi=e_0$.
Hence we infer that
$$
\varphi(S_1,\cdots,S_n)\psi(S_1,\cdots,S_n)=
\psi(S_1,\cdots,S_n)\varphi(S_1,\cdots,S_n)=I_{F^2(H_n)}.
$$
Therefore the operator $\varphi(S_1,\cdots,S_n)$ is invertible which
implies:
$$
\text{range\ }\varphi(S_1,\cdots,S_n)=F^2(H_n),\qquad\text{and}\leqno{(10)}
$$
$$
\varphi(S_1,\cdots,S_n)\quad\text{is bounded below.}\leqno{(11)}
$$

The relation (10) shows that $~\varphi~$ is outer and (11) implies that
there exists $\delta>0$ such that
$$
\|\varphi\otimes p\|_2\ge\delta\|p\|_2,
\quad\text{for any\ } p\in\P.\leqno{(12)}
$$

Conversely, assume $~\varphi\in F^\infty$ is outer such that (12)
holds. This implies that
$\varphi(S_1,\cdots,S_n)$ is invertible in $B(F^2(H_n))$. Let
$T\in B(F^2(H_n))$ be such that
$$
T\varphi(S_1,\cdots,S_n)=\varphi(S_1,\cdots,S_n)T=I_{F^2(H_n)}
$$
Then for any $~h\in F^2(H_n)~$ we have
$$
T(\varphi\otimes h)=\varphi\otimes T h=h.\leqno{(13)}
$$
Let $~h=e_0~$ and $~\psi=T e_0$. Then (13) gives us
$$
\varphi\otimes\psi=e_0\quad\text{and}\quad T(\varphi)=e_0.\leqno{(14)}
$$
We clearly have
$$
\varphi=e_0\otimes\varphi=(\varphi\otimes\psi)\otimes\varphi=\varphi\otimes
(\psi\otimes\varphi).
$$
Then according to (14) and (13) we have
$e_0=T(\varphi)=T(\varphi\otimes(\psi\otimes\varphi))=\psi\otimes\varphi$.
Therefore,
$$
\varphi\otimes\psi=\psi\otimes\varphi=e_0.
$$
We claim that $~\psi\in F^\infty$.
Let $~p\in\P$. It is clear that
$$
p=e_0\otimes p=(\varphi\otimes\psi)\otimes p=\varphi\otimes(\psi\otimes
p)
$$
Therefore, by (13)
$$
T(p)=T(\varphi\otimes(\psi\otimes p))=\psi\otimes p.
$$
Hence,
$$
\|\psi\|_\infty=\sup_{p\in(\P)_1}\|\psi\otimes p\|_2=
\sup_{p\in(\P)_1}\|T p\|_2=\|T\|<\infty.
$$
The proof is complete. $\qed$
\enddemo

\heading
3. Examples of inner and outer functions
\endheading
In this section we will present many examples of inner and outer
functions
in $F^\infty$. Recall that $\varphi\in F^\infty$ is inner iff
the map $h\to \varphi\otimes h$ is an isometry on $\F^2(\H_n)$,
(equivalently, $h\to h\otimes\tilde{\varphi}$ is an isometry
on $\F^2(\H_n)$), and $\varphi\in F^\infty$ is outer iff there
exists $p_n\in\F^2(\H_n)$ such that $\varphi\otimes p_n\to e_0$ in
$\F^2(H_n)$.

\medbreak

The first observation is that any inner (respectively outer) function
in $H^\infty$ becomes inner (respectively outer) in $F^\infty$ if
we substitute $z$ by $e_f$. We call these examples ``inherited''.

\proclaim{Example 1} Inherited inner and outer functions.
\endproclaim

For any $f\in\F$, $f\not=0$ let $V_f:H^2\to\F^2(H_n)$ be the isometry
defined by $V_f z^k=e_f^k$. We denote $V_f\varphi=\varphi_f$.

\proclaim{Proposition 3.1} $V_f:H^\infty\to F^\infty$ is an isometry.
Moreover,
if $\varphi\in H^\infty$ is inner then $\varphi_f$ is inner in
$F^\infty$,
and if $\varphi\in H^\infty$ is outer then $\varphi_f$ is outer
in $F^\infty$.
\endproclaim

\demo{\bf Proof} Let $\F_f^2=[e_0,e_f,e_f^2,\cdots]$. Any $h\in\F$
satisfies
$e_h=e_f^k\otimes e_g$ for some $k\in\N_0$ and $g\in\F$ that
``does-not-start'' from $f$ (i.e., $S_f^*e_g=0$, see (4) for the
notation
of $S_f$). We decompose $\F^2(\H_n)$ as
$$\F^2(\H_n)=\bigoplus_{{g\in\F\atop S_f^*e_g=0}}~
[\F_f^2\oplus e_g].\leqno{(15)}$$
Since $\varphi_f\in\F_f^2$ we have that for any $g\in\F$ that does not
start from $f$,
$$\varphi_f\otimes[\F_f^2\otimes e_g]\subset\F_f^2\otimes e_g.$$
Moreover, it is easy to see that
$$\sup_{p\in(\F_f^2\otimes e_g)_1}\|\varphi_f\otimes p\|_2
=\sup_{p\in(\F_f^2)_1}\|\varphi_f\otimes p\|_2
=\sup_{p\in(H^2)_1}\|\varphi p\|_2=\|\varphi\|_\infty.$$
Combining this with the fact that $\varphi_f$ acts diagonally in the
decomposition of $\F^2(\H_n)$ given by (15) we conclude that
$\|\varphi_f\|_\infty=\|\varphi\|_\infty$.
\smallbreak

Suppose that $\varphi\in H^\infty$ is outer. Then we can find a
sequence
of $p_n\in H^2$ such that $\varphi p_n\to 1$ in $H^2$. Since $V_f$ is
an isometry on $H^2$ we get that
$V_f(\varphi p_n)=\varphi_f \otimes V_f(p_n)\to e_0$ in $\F^2(\H_n)$.
By Proposition 1.4
we get that $\varphi_f$ is outer in $F^\infty$.
\smallbreak

Suppose now that $\varphi\in H^\infty$ is inner. We want to prove that
whenever $h_1, h_2\in\F$ satisfy $h_1\not= h_2$, we will have that
$\varphi_f\otimes e_{h_1}\perp \varphi_f\otimes e_{h_2}$.
\smallbreak

Consider first the case $e_{h_1}=e_f^k$, $e_{h_2}=e_f^\ell$ for
$k\not=\ell$.
Since $\varphi\in H^\infty$ is inner then $\varphi z^k\perp\varphi
z^\ell$.
Hence using that $V_f$ is an isometry we get that
$\varphi_f\otimes e_{h_1}\perp\varphi_f\otimes e_{h_2}$.

Suppose now that $e_{h_i}\in\F_f^2\otimes e_{g_i}$, $i=1,2$ for some
$g_1, g_2\in\F$ that do not start from $f$.

If $g_1\not=g_2$ then $\varphi_f\otimes e_{h_i}\in\F_f^2\otimes
e_{g_i}$ for
$i=1,2$. Since these two subspaces are orthogonal we conclude that
$\varphi_f\otimes e_{h_1}\perp\varphi_f\otimes e_{h_2}$.

If $g_1=g_2$
then $e_{h_1}=e_f^k\otimes e_g$ and $e_{h_2}=e_f^\ell\otimes e_g$ for
some $k\not=\ell$. Hence,
$$\langle \varphi_f\otimes e_{h_1},\varphi_f\otimes e_{h_2}\rangle
=\langle \varphi_f\otimes e_f^k,\varphi_f\otimes e_f^\ell\rangle=0.$$
The last equality follows from the first case.\qed
\enddemo

\proclaim{Example 2} (Invertible elements)
If $\psi\in F^\infty$ is invertible, then $\psi$
is outer.
\endproclaim

In particular, the following examples are invertible (and hence
outer):
\roster
\item"{(i)}" For $\varphi_i\in F^\infty$, $\|\varphi_i\|_\infty<1$
for $i\leq k$, let
$$\psi=(e_0-\varphi_1)\otimes (e_0-\varphi_2)\otimes\cdots\otimes
(e_0-\varphi_k).$$
\item"{(ii)}" For $\varphi\in F^\infty$, and
$\varphi^n=\varphi\otimes\cdots\otimes\varphi$ ($n$-times), let
$$\psi=\exp \varphi=\sum_{k=0}^\infty{\varphi^n\over n!}.$$
\item"{(iii)}" For $\lambda\in\C^n$, $\|\lambda\|_2<1$ let
$\psi=z_\lambda$, where the $z_\lambda$'s are defined in Example 8.
\endroster

\bigbreak

We will present examples of inner functions now. The first example
is the simplest.

\proclaim{Example 3} (The monomials) For every $f\in\F$, $e_f$ is
inner.
\endproclaim

Examples 4 and 6 below appear in [A].

\proclaim{Example 4} For every $k\in\N$ let
$X_k=\hbox{span}~\{e_f:f\in F(k,\Lambda)\}$ (i.e., the span of the
monomials
with $k$ letters). Then any $x\in X_k$, $\|x\|_2=1$ is inner.
\endproclaim

\demo{\bf Proof}
The main point is that if $f_1,f_2\in F(k,\Lambda)$
and $g,h\in\F$ then $e_{f_1}\otimes e_g=e_{f_2}\otimes e_h$ iff
$f_1=f_2$ and $g=h$.

We can easily see that if $x=\sum_{f\in F(k,\Lambda)} a_f e_f\in X_k$
and $g,h\in\F$, $g\not= h$ then $x\otimes e_g$ is orthogonal to
$x\otimes e_h$.
If we require that $\|x\|_2=1$, $x$ is inner.\qed
\enddemo

\proclaim{Example 5} Let $e_{f_1}, \cdots, e_{f_k}$ be monomials such
that
the first letter of all of them are different. (i.e., $k\not= \ell$
implies
that $f_k(1)\not= f_\ell(1)$). Then for any $\sum_{i\leq k}|a_i|^2=1$
we have
that $\sum_{i\leq k}a_ie_{f_i}$ is inner.
\endproclaim

\demo{\bf Proof}
The proof is similar to that of Example 4. Since
the $f_i$'s start from different letters, for any $g,h\in\F$,
$e_{f_i}\otimes e_g\perp e_{f_j}\otimes e_h$ if $i\not= j$.

We can easily see that if $x=\sum_{i\leq k} a_i e_{f_i}$ and $g,
h\in\F$,
$g\not=h$ then $x\otimes e_g$ is orthogonal to $x\otimes e_h$. If we
require that $\|x\|_2=1$, $x$ is inner.\qed
\enddemo

\proclaim{Example 6} Let $\psi\in\F^2(H_{n-1})$, $\|\psi\|_2=1$,
(i.e., $\psi$ does not have any $e_n$). Then $\psi\otimes e_n$ is
inner.
\endproclaim

\demo{\bf Proof}
Let $\F_{n-1}$ be the set of words on the letters $e_1,\cdots,
e_{n-1}$.
Let $f_1,f_2\in\F_{n-1}$ and $g_1,g_2\in\F$. The main point is that
$$e_{f_1}\otimes e_n\otimes e_{g_1}=e_{f_2}\otimes e_n\otimes e_{g_2}
\quad\text{ if and only if }\quad f_1=f_2\text{ and }g_1=g_2.$$
To see this notice that the first time that $e_n$ appears in
$e_{f_1}\otimes e_n\otimes e_{g_1}$ is right after $e_{f_1}$.
Since the same is true for $f_2$ the words must
agree before $e_n$ (i.e., $f_1=f_2$) and after (i.e., $g_1=g_2$).

It is easy to prove now that if $g_1,g_2\in\F$, $g_1\not= g_2$,
then $\psi\otimes e_n\otimes e_{g_1}\perp \psi\otimes e_n\otimes
e_{g_2}$.
The normalization condition guarantees that $\psi\otimes e_n$ is
inner.\qed
\enddemo

\proclaim{Example 7} (M\"obius functions) For any $f\in\F$, and
$\mu\in\C$ with $|\mu|<1$, we have that
$\varphi(f,\mu)=(e_f-\mu)\otimes(1-\bar{\mu}e_f)^{-1}$ is inner.
\endproclaim

These are particular cases for Example 1 above.
The products of these elements can be viewed as analogue to
the Blashke products. However,
$[\varphi(f,\mu)\otimes\F^2(\H_n)]^\perp$ is always infinite
codimensional. In fact, since $(1-\bar{\mu}e_f)^{-1}$ is invertible and
$(1-\bar{\mu}e_f)^{-1}$, $(e_f-\mu)$ commute, we have that
$h\in [\varphi(f,\mu)\otimes\F^2(\H_n)]^\perp$ iff for every
$\psi\in \F^2(\H_n)$,
$\langle(e_f-\mu)\otimes\psi,h\rangle=0.$ Then, looking
only at the basic elements we get that
$$
h\in [\varphi(f,\mu)\otimes\F^2(\H_n)]^\perp
\Longleftrightarrow \forall g\in\F,~
\langle e_f\otimes e_g,h\rangle=\mu\langle e_g,h\rangle.\leqno{(16)}
$$
It is easy to see that
$$
h_f(\mu)=\sum_{k=0}^\infty \mu^k e_f^k
$$
satisfies $(16)$ and thus belongs to
$[\varphi(f,\mu)\otimes\F^2(\H_n)]^\perp$.

One can also check that if $g\in\F$ ``does-not-start'' from $f$
(i.e., $S_f^*e_g=0$, see (4) for the notation of $S_f$), then
$h_f(\mu)\otimes e_g\in [\varphi(f,\mu)\otimes\F^2(\H_n)]^\perp$.
Moreover, the span of these elements is dense in
$[\varphi(f,\mu)\otimes\F^2(\H_n)]^\perp$. We leave the details out.
\smallbreak

We need the following lemma in the next section.

\proclaim{Lemma 3.2} Let $\Omega\subset \{\mu\in\C:|\mu|<1\}$ be a set
with an
accumulation point in $\{\mu\in\C:|\mu|<1\}$. Then
$\bigcap\{\varphi(f,\mu)\otimes\F^2(\H_n):f\in\F,\mu\in\Omega\}=\{0\}.$
\endproclaim

\demo{\bf Proof}
Let
$\psi\in
\bigcap\{\varphi(f,\mu)\otimes\F^2(\H_n):f\in\F,\mu\in\Omega\}$
and fix $f\in\F$.
For every $\mu\in\Omega$, we have that
$\langle\psi,h_f(\mu)\rangle=0$. Hence,
$$\langle\psi,e_0\rangle+\mu \langle\psi,e_f\rangle+\mu^2
\langle\psi, e_f^2\rangle+\mu^3 \langle\psi,e_f^3\rangle+\cdots=0
\quad\quad\hbox{for all}\quad \mu\in\Omega.$$
Since the map
$\mu\to\sum_{k=0}^\infty \mu^k \langle\psi,e_f^k\rangle$ is analytic
on $\mu$ and the zeros accumulate inside the disk, the map is zero.
Hence $\langle\psi,e_0\rangle=\langle\psi,e_f\rangle=0$.
Since $f$ is arbitrary, we conclude
that $\psi=0$.\qed
\enddemo

The following example is more technical and it is not used in the
rest of the paper. However, it gives rise to non-trivial inner
functions and to an interesting ideal in $F^\infty$ that seems to
capture the ``non-commutativity'' of the product.

\proclaim{Example 8} 1-codimensional invariant subspaces.
\endproclaim

Let $\M$ be a 1-codimensional subpace of $\F^2(\H_n)$ invariant for
$S_1,\cdots,S_n$. Then $\M=[z]^\perp$ for some $z\in\F^2(\H_n)$ and
$[z]$ is invariant for $S_1^*,\cdots, S_n^*$. That is, for every
$i\leq n$ there exists $\lambda_i\in\C$ such that
$$S_i^*z=\lambda_i z.$$
Assume that $\langle z,e_0\rangle=1$.

For $f\in F(k,\Lambda)$ let
$\lambda_f=\lambda_{f(1)}\lambda_{f(2)}\cdots\lambda_{f(k)}$ and
$\lambda_0=1$. We claim that for any $f\in\F$
$$\langle z,e_f\rangle=\lambda_f\quad\hbox{ and then }\quad
z=\sum_{f\in\F}\lambda_fe_f.\leqno{(17)}$$
To see this notice that for every $i\leq n$,
$\langle z,e_i\rangle=\langle z,S_ie_0\rangle=\langle S_i^*z,e_0\rangle
=\lambda_i$. Similarly,
$\langle z,e_i\otimes e_j\rangle=\langle z,S_iS_je_0\rangle
=\langle S_j^*S_i^*z,e_0\rangle=\lambda_i\lambda_j$. Proceeding
inductively we get (17). Since the $\lambda_i$'s determine $z$ and $\M$
we will denote them by $z_\lambda$ and $\M_\lambda$ from now on.

It is not hard to see that
$$z_\lambda=\sum_{f\in\F}\lambda_fe_f=1+\sum_{k=1}^\infty
(\lambda_1e_1+\lambda_2e_2\cdots+\lambda_ne_n)^k.$$
Then a necessary and sufficient condition for $z_\lambda\in\F^2(\H_n)$
is that $\|\lambda\|_2=\sqrt{\sum_{i\leq n}|\lambda_i|^2}<1$. Moreover,
by Example 4 and Proposition 1.6,
$\|\lambda_1e_1+\cdots\lambda_ne_n\|_2
=\|\lambda_1e_1+\cdots\lambda_ne_n\|_\infty$ and then $z_\lambda\in
F^\infty$.

We have thus proved:

\proclaim{Proposition 3.3} $\M$ is a 1-codimensional invariant subspace
for
$S_1,\cdots,S_n$ if and only if $\M=[z_\lambda]^\perp$ for some
$\lambda\in\C^n$, $\|\lambda\|_2<1$.
\endproclaim

It follows from Proposition 2.3 that
$\M_\lambda=\bigoplus_{j\in J}[\F^2(\H_n)\otimes\tilde{\varphi}_j]$ for
some orthogonal inner functions $\varphi_j$, $j\in J$.
Recall that $|J|$, the cardinality of $J$, is the dimension
of the wandering subspace
$\L_\lambda=\M_\lambda\ominus
[S_1\M_\lambda\oplus\cdots\oplus S_n\M_\lambda]$ for
$S_1|_{\M_\lambda}, \cdots, S_n|_{\M_\lambda}$. Moreover, whenever
$\tilde{\varphi}\in\L_\lambda$ satisfies $\|\tilde{\varphi}\|_2=1$ we
have
that $\varphi$ is inner.

Let $Q_\lambda$ be the orthogonal projection onto $\M_\lambda$ and
$P_\lambda$ the orthogonal projection onto $\L_\lambda$. One can
easily check that
$P_\lambda=Q_\lambda-\sum_{i\leq n}S_iQ_\lambda S_i^*$ and that
$$\eqalign{
Q_\lambda e_f&=e_f-{\lambda_f\over\|z_\lambda\|_2^2}z_\lambda\quad
\hbox{ for any }\quad f\in\F,\cr
P_\lambda e_i&={1\over\|z_\lambda\|_2^2}(e_i-\lambda_i)\otimes
z_\lambda
\quad\hbox{ for }\quad i\leq n,\cr
P_\lambda e_0&=e_0-{1\over\|z_\lambda\|_2^2}z_\lambda.\cr}$$
If $e_f=e_i\otimes e_g$ for some $i\leq n$ and $g\in\F$ we can also
check that
$P_\lambda e_f={\lambda_g\over\|z_\lambda\|_2^2}
(e_i-\lambda_i)\otimes z_\lambda$.

It is easy to see that (17) implies that $\tilde{z}_\lambda=z_\lambda$,
since $\lambda_f=\lambda_{\tilde{f}}$. Then we get

\proclaim{Proposition 3.4}
$\L_\lambda=\hbox{span}\,\{\tilde{\varphi}_0,
\tilde{\varphi}_1,\cdots,\tilde{\varphi}_n\}$ where
$$\varphi_0=a_0\bigl[e_0-{1\over\|z_\lambda\|_2^2}z_\lambda\bigr],\quad
\hbox{ and }\quad
\varphi_i=a_i\bigl[z_\lambda\otimes (e_i-\lambda_i)\bigr]\quad
i=1,2\cdots,n.$$
The $a_i\in\C$ are chosen so that $\|\varphi_i\|_2=1$ for
$i=0,1,\cdots, n$.
Moreover, the $\varphi_i$'s and any normalized linear combination of
them
are inner function.
\endproclaim

\noindent{\bf Remarks.} (1) We have that
$\M_\lambda=\F^2(\H_n)\otimes\tilde{\varphi}_0+
\F^2(\H_n)\otimes\tilde{\varphi}_1+\cdots+
\F^2(\H_n)\otimes\tilde{\varphi}_n$. However, the $\varphi_i$'s are not
orthogonal inner functions.
\smallbreak
(2) We also have that $\psi\in\M_\lambda$ if and only if
$\langle\psi,z_\lambda\rangle=0$ if and only if
there exist $g_i\in\F^2(\H_n)$, $i=0,1,\cdots,n$ such that
$$\psi=g_0\otimes\tilde{\varphi}_0+g_1\otimes\tilde{\varphi}_1+
\cdots+g_n\otimes\tilde{\varphi}_n.$$
This can be view as a ``weak-factorization'' result.
\smallbreak
(3) It is easy to compute that $\|z_\lambda\|_2^2={1\over
1-\|\lambda\|_2^2}$.

\medbreak
If $\M$ is an invariant subspace for $S_1,\cdots, S_n$ and $p\in\P$
we always have (Proposition 2.3) that $p\otimes\M\subset\M$. In
general we do not have that $\M\otimes p\subset \M$. However, this is
always
true for the 1-codimensional invariant subspaces.

\proclaim{Lemma 3.5} Let $\M_\lambda=[z_\lambda]^\perp$ for some
$\lambda\in\C^n$, $\|\lambda\|_2<1$. For any $p\in\P$ we have that
$p\otimes\M_\lambda\subset\M_\lambda$ and $\M_\lambda\otimes p\subset
\M_\lambda$.
\endproclaim

\demo{\bf Proof}
The first inclusion is clear. If $\varphi\in\M_\lambda$ and
$p\in\P$ then $\langle p\otimes\varphi,z_\lambda\rangle=0$.

The key point for the other one is that $\varphi\in\M_\lambda$ if and
only if
$\tilde{\varphi}\in\M_\lambda$. To see this recall that
$z_\lambda=\tilde{z}_\lambda$. Then
$$\langle \tilde{\varphi}, z_\lambda\rangle
=\langle \varphi,\tilde{z}_\lambda\rangle
=\langle \varphi, z_\lambda\rangle.$$
If $p\in\P$ and $\varphi\in\M_\lambda$ we have that
$$\langle \varphi\otimes p,z_\lambda\rangle
=\langle \varphi\otimes p,\tilde{z}_\lambda\rangle
=\langle \widetilde{\varphi\otimes p},z_\lambda\rangle
=\langle \tilde{p}\otimes\tilde{\varphi},z_\lambda\rangle=0.\qed$$
\enddemo

Let $\Nu=\bigcap_{\|\lambda\|_2<1}\M_\lambda$ and
$\Nu^\infty=\Nu\bigcap F^\infty$. $\Nu$ is a ``2-sided'' invariant
subspace and $\Nu^\infty$ is a two-sided ideal in $F^\infty$.

There are plenty of elements in $\Nu$. For example
$\varphi=e_1\otimes e_2-e_2\otimes e_1\in\Nu$. For every
$\lambda\in\C^n$,
$\|\lambda\|_2<1$ we have
$$\langle \varphi,z_\lambda\rangle
=\langle e_1\otimes e_2,z_\lambda\rangle-
\langle e_2\otimes e_1,z_\lambda\rangle
=\lambda_1\lambda_2-\lambda_2\lambda_1=0.$$
More generally if $f\in F(k,\Lambda)$ and $\pi\in\Pi_k$, ($\pi$ is a
permutation on $\{1,2,\cdots,k\}$), let $\pi(f)\in F(k,\Lambda)$ be
defined by $\pi(f)(i)=f(\pi(i))$. Then $e_f-e_{\pi(f)}\in\Nu$.

\proclaim{Proposition 3.6} $F^\infty/\Nu^\infty$ is a commutative
algebra.
\endproclaim

\demo{\bf Proof} Let $\varphi=\sum_{f\in\F}a_f e_f$,
$\psi=\sum_{g\in\F}b_g e_g$ be elements in $F^\infty$. We want
to prove that $\varphi\otimes\psi-\psi\otimes\varphi\in\Nu$.
Notice that
$$\varphi\otimes\psi-\psi\otimes\varphi=
\sum_{f\in\F}\sum_{g\in\F}a_fb_g[e_f\otimes e_g-e_g\otimes e_f].$$
Let $\lambda\in\C^n$, $\|\lambda\|_2<1$. Since
$\langle e_f\otimes e_g-e_g\otimes e_f, z_\lambda\rangle=0$
we get that
$\langle \varphi\otimes\psi-\psi\otimes\varphi,z_\lambda\rangle=0.$
Since $\lambda$ is arbitrary we finish.\qed
\enddemo

\heading
4. Reflexivity results
\endheading

Let $~\H~$ be a Hilbert space and $~B(\H)~$ be the algebra of all
bounded
operators on $~\H$. If $~A\in B(\H)~$ then the set of all invariant
subspaces
of $~A~$ is denoted by $\text{Lat~}A$. For any $~\U\subset B(\H)~$ we
define
$$
\text{Lat~}\U=\bigcap_{A\in\U}\text{Lat~}A.
$$

If $~\S~$ is any collection of subspaces of $~\H$, then
$\text{Alg~}\S:=\{A\in B(\H):\ \S\subset\text{Lat~}A\}.$
The algebra $~\U\subset B(\H)~$ is reflexive if $\U=\text{Alg Lat~}\U.$

The main theorems of this section are the following.

\proclaim{Theorem 4.1} The algebra
$\{\varphi(S_1,S_2,\cdots,S_n):\varphi\in F^\infty\}$ is reflexive.
\endproclaim

\proclaim{Theorem 4.2} If $~\U~$ is a strongly closed subalgebra of
$~\{\varphi(S_1,\cdots,S_n);\  \varphi\in F^\infty\}~$
containing the identity then $~\U~$ is reflexive.
\endproclaim

\demo{\bf Proof of Theorem 4.1}
Let $A\in \hbox{Alg Lat}\, F^\infty$.
For every $\varphi$ inner we have that
$F^2(H_n)\otimes\tilde{\varphi}\in\hbox{Lat}\,F^\infty$. Then
$A [F^2(H_n)\otimes\tilde{\varphi}]\subset
F^2(H_n)\otimes\tilde{\varphi}$.
In particular we have
$$A\tilde{\varphi}=\psi_\varphi\otimes\tilde{\varphi}\quad\quad
\hbox{for some}\quad \psi_\varphi\in F^2(H_n).\leqno{(18)}$$

\smallbreak
For every $f\in\F$ we can find $\psi_f$ such that
$$A e_f=\psi_f\otimes e_f.$$
(Notice that both $e_f$ and $\tilde{e}_f$ are inner and we do not have
to carry the tilde).

\smallbreak
Let $k\geq 1$. Example 4 of Section 3 tells that
$x_k=a_k\sum_{f\in F(k,\Lambda)} e_f$, where
$a_k=[\hbox{card}(F(k,\Lambda))]^{-1/2}$, is inner
(notice that $x_k=\tilde{x}_k$).
Using $(18)$ we have that
$$Ax_k=\psi_k\otimes x_k=a_k\sum_{f\in F(k,\Lambda)}\psi_k\otimes e_f.
\leqno{(19)}$$
On the other hand, we have that
$$Ax_k=A\biggl(a_k\sum_{f\in F(k,\Lambda)} e_f\biggr)=
a_k\sum_{f\in F(k,\Lambda)} Ae_f=a_k\sum_{f\in
F(k,\Lambda)}\psi_f\otimes e_f.
\leqno{(20)}$$
We claim that for every $f\in F(k,\Lambda)$, $\psi_k=\psi_f$. Fix
$f\in F(k,\Lambda)$ and let $P_f$ be the orthogonal projection onto
$\F^2(\H_n)\otimes e_f$. From (19) we get that
$P_f(Ax_k)=a_k\psi_k\otimes e_f$ (the other terms are zero), and from
(20) we get that $P_f(Ax_k)=a_k\psi_f\otimes e_f$. Therefore,
$\psi_k\otimes e_f=\psi_f\otimes e_k$ and then $\psi_k=\psi_f$.

\smallbreak
Let $k>1$. Example 5 of Section 3
tells us that $y={1\over\sqrt{2}}[e_1+e_2^k]$
is inner (notice that $y=\tilde{y}$). Using $(18)$ we have that
$$Ay=\psi_y\otimes y={1\over\sqrt{2}}[\psi_y\otimes e_1+\psi_y\otimes
e_2^k].
\leqno{(21)}$$
On the other hand,
$$Ay={1\over\sqrt{2}}[Ae_1+Ae_2^k]={1\over\sqrt{2}}
[\psi_1\otimes e_1+\psi_k\otimes e_2^k].\leqno{(22)}$$
The last equality is clear since $e_1$ has {\sl one} letter and
$e_2^k$ has $k$ letters. Combining (21) and (22) we conclude that
$\psi_k=\psi_1$.

\smallbreak
Summarizing we have that $Ae_0=\psi_0$ for some $\psi_0\in F^2(H_n)$,
and
if $f\in F(k,\Lambda)$, $k\geq 1$, then $Ae_f=\psi_1\otimes e_f$. It is
easy to see that $\psi_1\in F^\infty$.

Let $B=A-\psi_1(S_1,S_2,\cdots, S_n)$. We still have that
$B\in\hbox{Alg Lat}\,F^\infty$ and if we set $\psi=\psi_0-\psi_1$,
$$\eqalign{Be_0&=\psi\cr
	   Be_f&=0\quad\quad\hbox{if}\quad f\not=0.\cr}$$
We want to prove that $B=0$. This gives us that $A=\psi_1\in
F^\infty$.

\smallbreak
Equation $(18)$ applies to $B$ as well. That is, if $\varphi$ is inner,
then $B\tilde{\varphi}=\psi_\varphi\otimes\tilde{\varphi}$. On the
other
hand it is clear that
$B\tilde{\varphi}=\langle e_0,\tilde{\varphi}\rangle\psi$.
Hence, if $\langle e_0,\tilde{\varphi}\rangle\not=0$, we have that
$$
\psi={1\over \langle e_0,\tilde{\varphi}\rangle }
			       \psi_\varphi\otimes \tilde{\varphi}
\in F^2(H_n)\otimes\tilde{\varphi},\quad\hbox{and}\quad
\tilde{\psi}\in\varphi\otimes\F^2(\H_n).
$$

It is easy to see that $\varphi(f,\mu)$ (the M\"obius maps of Example
7,
Section 3) satisfy $\langle e_0,\tilde{\varphi}(f,\mu)\rangle\not=0$
for any
$f\in\F$ and $0<|\mu|<1$. This implies that
$$
\tilde{\psi}\in\bigcap\{\varphi(f,\mu)\otimes \F^2(\H_n): f\in\F,
0<|\mu|<1\}.
$$
By Lemma 3.2 we conclude that $\psi=0$. \qed
\enddemo

\noindent{\bf Remark} The previous proof does not work in the
commutative
case. We are using the fact that $e_1$ and $e_2$ are non-commutative
to conclude that ${1\over\sqrt{2}}[e_1+e_2^k]$ is inner. However,
${1\over\sqrt{2}}[z+z^k]$ is not inner in $H^\infty$.

The following notation will be useful in the proof of Theorem 4.2:
If $~A\in B(\H)~$ and
$~m~$ is a positive integer, then $~\H^{(m)}~$ denotes the direct sum
of
$~m~$ copies of $~\H~$ and $~A^{(m)}~$ stands for the direct sum of
$~m~$ copies of $~A$.
If $~\U\subset B(\H)$, then $~\U^{(m)}:=\{A^{(m)}:\ A\in\U\}$.

According to [RR, Theorem 7.1], if $~\U~$ is an algebra of operators
containing the identity, then the closure of $~\U~$ in the strong
operator topology is
$$
\{B:\text{Lat~}\U^{(m)}\subset\text{Lat~}B^{(m)}\quad\text{for\ }
m=1,2,3,\cdots\}.
$$

To prove Theorem 4.2 we need the following

\proclaim{Theorem 4.3}
If $~A, B\in \{\varphi(S_1,\cdots,S_n);\quad\varphi\in F^\infty\}~$
are such that
$\text{Lat~}A\subset\text{Lat~}B$, then
$$
\text{Lat~}A^{(m)}\subset \text{Lat~}B^{(m)},\quad\text{for any\ }
m=1,2,\cdots.
$$
\endproclaim

\demo{\bf Proof}
It will enough to show that every cyclic invariant subspace of
$~A^{(m)}~$ is
invariant under $~B^{(m)}$. Let  $~x\not=0~$ be an element in
$~F^2(H_n)^{(m)}~$ and let $~\M=\bigvee\limits_{f\in\F}S_f^{(m)}x$.
According to Lemma 2.5 there is a unitary operator
$~U:F^2(H_n)\to\M~$ such that
$$
S_i^{(m)}|_{\M}=U S_i U^{-1},\quad\text{for any\ }i=1,2,\cdots.
$$
Since $~A~$ and $~B~$ are strong limits of polynomials in
$~S_1,\cdots,S_n~$ it follows that
$$\eqalign{A^{(m)}|_{\M}&=U A U^{-1},\quad\hbox{ and}\cr
	   B^{(m)}|_{\M}&=U B U^{-1}.\cr}\leqno{(23)}$$

Let $~\G~$ be the smallest invariant subspace of $~A^{(m)}~$ containing
$x$, i.e.,
$$
\G=\bigvee\limits_{f\in\F}A_f^{(m)}x.
$$
Since $~A^{(m)}=\varphi(S_1^{(m)},\cdots,S_n^{(m)})~$
for some $~\varphi\in F^\infty~$ it follows that $~\G\subset \M$.
The relation (23) and the fact that $~A^{(m)} \G\subset\G~$
implies $~U A U^{-1}(\G)\subset\G~$ whence
$~A(U^{-1}\G)\subset U^{-1}\G$.
But $~\text{Lat~}A\subset\text{Lat~}B~$ implies
$~B(U^{-1}\G)\subset U^{-1}\G$.
Therefore, $U B U^{-1}(\G)\subset\G~$ which together with (23) show
that
$~B^{(m)}\G\subset\G$. \qed
\enddemo

\proclaim{Corollary 4.4}
Let $~A, B\in\{\varphi(S_1,\cdots,S_n):\quad\varphi\in F^\infty\}~$
be such that
$~\text{Lat~}A\subset\text{~Lat~}B$, then $~B~$ belongs to the strongly
closed algebra
generated by $~A~$ and the identity.
\endproclaim

The proof of the following theorem is analogue to the proof of the
Theorem 4.3
if we replace $~A~$ with any subset
$\A\subset\{\varphi(S_1,\cdots,S_n):\quad\varphi\in F^\infty\}$.
We omit the proof.

\proclaim{Theorem 4.5}
If $~B\in\{\varphi(S_1,\cdots,S_n):\quad\varphi\in F^\infty\}~$
is such that $~\text{Lat~}\A\subset\text{~Lat~}B~$ then
$$
\text{Lat~}\A^{(m)}\subset\text{~Lat~}B^{(m)},\quad\text{for
any~}m=1,2,\cdots.
$$
Moreover, $~B~$ belongs to the strongly closed algebra generated by
$~\A~$ and the identity.
\endproclaim

\demo{\bf Proof of Theorem 4.1}
Let $~B\in B(\F^2(\H_n))~$ be such that
$$
\text{Lat~}\U\subset\text~{Lat~}B.
$$
It is clear that
$$
\text{Lat~}\{S_1,\cdots,S_n\}\subset\text{~Lat~}\{\varphi(S_1,\cdots,S_n);\quad
\varphi\in F^\infty\}\subset\text{~Lat~}\U\subset\text{~Lat~} B.
$$
According to Theorem 4.1, we have
$$
B\in\{\varphi(S_1,\cdots,S_n);\quad\varphi\in F^\infty\}.
$$
Now Theorem 4.5 implies that $B\in\U$. \qed
\enddemo

Let us recall from [Po2] that an operator $T\in B(\F^2(\H_n))$
is called multi-analytic if $TS_i=S_iT$ for each
$i\in\Lambda=\{1,2,\cdots,n\}$. The following result is an easy
consequence of Theorem 4.2 and the characterization of the
multi-analytic
operators in terms of their symbols [Po4].

\proclaim{Corollary 4.6} Any strongly closed subalgebra of
multi-analytic operators containing the identity is reflexive.
\endproclaim

Let us remark that in the particular case when $\Lambda=\{1\}$ we find
again
the result of Sarason [S].

\heading
5. Open Questions
\endheading

We finish this paper with some questions.

Let $\hbox{Inv}~(F^\infty)$ be the group of
invertible elements in $F^\infty$, and $\G_0$ be the connected
component in  $\hbox{Inv}~(F^\infty)$ which contains the identity.
From the general theory of Banach algebras [D] we know that the
collection
of finite products of elements in $\exp\F^\infty$ is $\G_0$ and that
$\hbox{Inv}~(F^\infty)/\G_0$ is a discrete group.

\proclaim{Problem 1} Characterize $\hbox{Inv}~(F^\infty)/\G_0$.
\endproclaim

\proclaim{Problem 2} Characterize the proper maximal invariant
subspaces
for the $\Lambda$-orthogonal shift of $\F^2(\H_n)$.
\endproclaim

It is clear that all of the $\M_\lambda$ of Example 8, Section 3 are
maximal.
However, there are many invariant subspaces for the shift that are not
inside any of them.

Take, for instance, $e_f=e_1\otimes e_2\otimes\cdots\otimes e_n$
and $\mu\in\C$, ${1\over\sqrt{n}}<|\mu|<1$. We claim that
$\M=\F^2(\H_n)\otimes\tilde{\varphi}(f,\mu)$ is not a subset
of $\M_\lambda$ for any $\lambda\in\C^n$, $\|\lambda\|_2<1$. (See
Example 6, Section 3 for the notation).

It is clear that $\M\subset\M_\lambda$ if and only if
$z_\lambda\in\M^\perp$ and since $z_\lambda=\tilde{z}_\lambda$ it is
equivalent to $z_\lambda\in[\varphi(f,\mu)\otimes\F^2(\H_n)]^\perp.$
Using (16) we get that $\M\subset\M_\lambda$ if and only if
$$\eqalign{
\forall g\in\F,\quad&
	    \langle e_f\otimes e_g\rangle=\mu\langle
	    e_g,z_\lambda\rangle,
		    \quad\hbox{ equivalently }\cr
\forall g\in\F,\quad & \lambda_f\lambda_g=\mu\lambda_g.\cr}$$
Since $\lambda_0=1$, we have that
$\lambda_f=\lambda_1\lambda_2\cdots\lambda_n=\mu$. However,
$\sum_{i\leq n}|\lambda_i|^n<1$ implies that
$$|\mu|=|\lambda_1\lambda_2\cdots\lambda_n|<{1\over\sqrt{n}},$$
contradicting the condition on $\mu$.

\proclaim{Problem 3} Characterize the proper maximal subspaces of
$\F^2(\H_n)$ of the form $\varphi\otimes\F^2(\H_n)$, $\varphi$ inner.
\endproclaim

It follows from Corollary 2.7 that if $\varphi_1,~\varphi_2$ are inner
and $\varphi_1\otimes\F^2(\H_n)\subset\varphi_2\otimes\F^2(\H_n)$ then
there exists $\varphi_3$ inner such that
$\varphi_1=\varphi_2\otimes\varphi_3$. Then the maximal subspaces
of the form $\varphi\otimes\F^2(\H_n)$ correspond to ``prime'' inner
functions.

\proclaim{Problem 4} Is every inner function $\varphi$ a product of
``prime'' ones? The question referes mainly to convergence.
\endproclaim

\enddocument

\enddocument